\documentclass{article}

\usepackage{amssymb}

\title{\bf Some developments in vertex operator algebra theory, old
and new}

\author{James Lepowsky\footnote{The author is partially supported by
NSF grant DMS--0401302.}}

\date{}

\begin{document}

\bibliographystyle{alpha}
\maketitle

\newtheorem{theorem}{Theorem}

\newtheorem{definition}
   {Definition}

\begin{abstract}
In this exposition, I discuss several developments in the theory of
vertex operator algebras, and I include motivation for the definition.
\end{abstract}

This is a more detailed version of the talk I gave at the Conference
on Lie Algebras, Vertex Operator Algebras and Their Applications.  I
thank Yi-Zhi Huang and Kailash Misra very much for organizing this
conference.

I would like to motivate the concept of vertex (operator) algebra,
including the definition; to discuss some of the main sources of the
theory, including Lie algebras and partition identities, the
``monstrous moonshine'' problem, and string theory and conformal field
theory; and to mention a selection of developments.  (This writeup is
not intended to be a comprehensive survey.)  Some of this talk is
drawn from the introductory material in \cite{FLM3}--\cite{FLM5},
\cite{DL2} and \cite{LL}, as well as from earlier expositions such as
\cite{LW5}, \cite{L1}, \cite{L2}, \cite{HL7}, \cite{HL8} and
\cite{L4}.

The mathematical notion of ``vertex algebra'' was introduced by
R. Borcherds \cite{B1}.  A variant of it, a notion of ``vertex
operator algebra,'' was introduced in \cite{FLM5}.  These notions are
algebraic formulations of concepts that had been developed by many
string theorists, conformal field theorists and quantum field
theorists, and formalized in \cite{BPZ} as certain ``operator
algebras,'' later called ``chiral algebras'' in physics.

Vertex (operator) algebra theory is inherently ``nonclassical,'' in
the same spirit in which string theory in physics is nonclassical and
also in the same spirit in which the sporadic finite simple groups in
mathematics are nonclassical.  String theory, its initial version
having been introduced in the late 1960s, is based on the premise that
elementary particles manifest themselves as ``vibrational modes'' of
fundamental strings, rather than points, moving through space
according to quantum-field-theoretic principles.  A string sweeps out
a two-dimensional ``world-sheet'' in space-time, and it is fruitful to
focus on the case in which the surface is a Riemann surface locally
parametrized by a complex coordinate.  The resulting two-dimensional
conformal (quantum) field theory has been studied extensively.

Meanwhile, in mathematics, the vast program to classify the finite
simple groups (see \cite{Gor}) had reached a dramatic point around
1980.  What turned out to be the largest finite simple group that is
sporadic (i.e., not belonging to one of the infinite families), the
Fischer-Griess Monster $\mathbb{M}$, had been predicted by B.~Fischer
and R.~Griess and was constructed by Griess \cite{G} as a symmetry
group (of order about ${10}^{54}$) of a remarkable new commutative but
very, very highly nonassociative, seemingly ad-hoc, algebra $\mathbb{B}$
of dimension 196,883.  The ``structure constants'' of the Griess
algebra $\mathbb{B}$ were ``forced'' by expected properties of the
conjectured-to-exist Monster.  It was proved by J.~Tits that 
$\mathbb{M}$ is actually the full symmetry group of $\mathbb{B}$.

A bit earlier (1978--79), J.~McKay, J.~Thompson, J.~Conway and
S.~Norton (see especially Conway-Norton \cite{CN}) had discovered
astounding ``numerology'' culminating in the ``monstrous moonshine''
conjectures relating the not-yet-proved-to-exist Monster $\mathbb{M}$ to
modular functions in number theory, namely:

There should exist a (natural) infinite-dimensional $\mathbb{Z}$-graded
module for $\mathbb{M}$ (i.e., representation of $\mathbb{M}$)
$$V=\bigoplus_{n=-1,0,1,2,3,\dots}V_n$$
such that
\begin{equation}\label{gradeddimV=J}
\sum_{n=-1,0,1,2,3,\dots}({\rm dim}\;V_n)q^n=J(q),
\end{equation}
where
\begin{equation}\label{J}
J(q)=q^{-1}+0+196884q+{\mbox {higher-order terms}},
\end{equation}
the classical modular function with its constant term set to 0.
$J(q)$ is the suitably normalized generator of the field of
$SL(2,\mathbb{Z})$-modular invariant functions on the upper half-plane,
with $q=e^{2\pi i \tau}$, $\tau$ in the upper half-plane.  (Note:
196884=196883+1.)  The existence of such a structure $V$ was
conjectured by McKay and Thompson.

More generally, as Conway and Norton conjectured, for every $g \in
\mathbb{M}$ (not just $g=1$), the generating function
\begin{equation}\label{gradedtraceV}
\sum_{n=-1,0,1,2,3,\dots}({\rm tr}\;g|_{V_n})q^n
\end{equation}
should be the analogous ``Hauptmodul'' for a suitable discrete
subgroup of $SL(2,\mathbb{R})$, a subgroup having a fundamental
``genus-zero property,'' so that its associated field of
modular-invariant functions has a single generator (a Hauptmodul).
(The left-hand side of (\ref{gradeddimV=J}) is the {\it graded
dimension} of the graded vector space $V$, and (\ref{gradedtraceV}) is
the {\it graded trace} of the action of $g$ on the graded space $V$;
the graded dimension is of course the graded trace of the identity
element $g=1$.)  The Conway-Norton conjecture subsumed a remarkable
coincidence that had been noticed earlier---that the 15 primes giving
rise to the genus-zero property (see A.~Ogg \cite{O}) are precisely
the primes dividing the order of the (conjectured-to-exist) Monster.

Proving these conjectures would give a remarkable connection between
classical number theory and ``nonclassical'' sporadic group theory.
The existence of a structure $V$ was soon essentially (but
nonconstructively) proved by Thompson, A.~O.~L.~Atkin, P.~Fong and
S.~Smith.  After Griess constructed $\mathbb{M}$, with I.~Frenkel and
A.~Meurman \cite{FLM2} we (constructively) proved the McKay-Thompson
conjecture, that there should exist a natural (whatever that was going
to mean) infinite-dimensional $\mathbb{Z}$-graded $\mathbb{M}$-module $V$
whose graded dimension is $J(q)$, as in (\ref{gradeddimV=J}).  The
graded traces of some, but not all, of the elements of the
Monster---the elements of an important subgroup of $\mathbb{M}$, namely,
a certain involution centralizer involving the largest Conway sporadic
group $Co_1$---were consequences of the construction, and these graded
traces were indeed (suitably) modular functions \cite{FLM2}.  We
called this $V$ ``the moonshine module $V^{\natural}$'' because of its
naturality (although the contruction of $V^{\natural}$ is not short).
See J. Tits (\cite{Ti1}, \cite{Ti2}) for discussions of the
construction of the Monster and of the moonshine module.

The construction \cite{FLM2} heavily used a number of different types
of ``vertex operators,'' all of them recently constructed, or newly
constructed as steps in \cite{FLM2}, along with their algebraic
structure and relations.  These were needed for the construction of
the structure $V^{\natural}$ itself, including a natural ``algebra of
vertex operators'' acting on it.  They were also needed for the
construction of a natural infinite-dimensional ``affinization'' of the
Griess algebra $\mathbb{B}$ acting on $V^{\natural}$.  This
``affinization,'' which was part of the new algebra of vertex
operators, is analogous to, but more subtle than, the notion of affine
Lie algebra, an example of which is discussed below.  More precisely,
the vertex operators were needed for a ``commutative affinization'' of
a certain natural 196884-dimensional enlargement ${\mathcal B}$ of
$\mathbb{B}$, with an identity element (rather than a ``zero'' element)
adjoined to $\mathbb{B}$.  This enlargement ${\mathcal B}$ naturally
incorporated the Virasoro algebra---the central extension of the Lie
algebra of formal vector fields on the circle---acting on
$V^{\natural}$; for us, the Virasoro algebra arose not because of its
role as a fundamental ``symmetry algebra'' in string theory but rather
because of the fact that a natural identity element for the algebra
was ``forced'' on us by our construction.  The vertex operators were
also needed for a natural ``lifting'' of Griess's action of $\mathbb{M}$
{}from the finite-dimensional space $\mathbb{B}$ to the
infinite-dimensional structure $V^{\natural}$, including its algebra
of vertex operators and its copy of the affinization of ${\mathcal
B}$.  Thus the Monster was now realized as the symmetry group of a
certain explicit ``algebra of vertex operators'' based on an
infinite-dimensional $\mathbb{Z}$-graded structure whose graded
dimension is the modular function $J(q)$.

Griess's construction of $\mathbb{B}$ and of $\mathbb{M}$ acting on 
$\mathbb{B}$ was a crucial guide for us, 
although we did not start by using his
construction; rather, we recovered it, as a finite-dimensional
``slice'' of a new infinite-dimensional construction using vertex
operator considerations.  In fact, our presentation of the Griess
algebra (see \cite{FLM1}, \cite{FLM2}, \cite{FLM5}) was short and
entirely canonical, and involved no choices or guesses of signs or
structure constants; one ``reads'' this algebra as in \cite{FLM1} from
canonical untwisted and twisted vertex operator structures newly
constructed starting from the Leech lattice (mentioned below).  As
presented in \cite{FLM1}, the 196884-dimensional algebra is simply the
direct sum of the weight-two subspace of the canonical
involution-fixed subspace of the untwisted Leech lattice vertex
operator structure with the analogous subspace of a canonically
twisted vertex operator structure; the commutative nonassociative
algebra structure and natural ``associative'' symmetric bilinear form
structure on ${\mathcal B}$ are essentially described ``in words.''
The initally strange-seeming {\it finite-dimensional} Griess algebra
was now embedded in a natural new {\it infinite-dimensional} space on
which a certain algebra of vertex operators acts, via a new kind of
``generalized commutation relation'' (relations of this type are
discussed below); such relations are what gave the commutative
affinization mentioned above.  At the same time, the Monster, a {\it
finite} group, took on a new appearance by now being understood in
terms of a natural {\it infinite-dimensional} structure.  The
very-highly-nonassociative Griess algebra, or rather, from our
viewpoint, the natural modification of the Griess algebra, with an
identity element adjoined, coming from a ``forced'' copy the Virasoro
algebra, became simply the conformal-weight-two subspace of an algebra
of vertex operators of a certain ``shape.''  The word ``simply''
refers to the ease of defining a commutative nonassociative algebra
with an associative symmetric bilinear form (generalizing the Griess
algebra with identity element adjoined, for the special case of
$V^{\natural}$) in the new general context of algebras of vertex
operators of ``shape'' similar to that of $V^{\natural}$ (as was
explained in \cite{FLM2} and \cite{FLM5}); the actual construction of
the particular algebra $V^{\natural}$ remains complex.  In any case,
the largest sporadic {\it finite} simple group, the Monster, was
``really'' {\it infinite-dimensional}.

In the expansion (\ref{J}), the constant term of $J(q)$ is zero, and
this choice of constant term, which is not uniquely determined by
number-theoretic principles, is not traditional in number theory.  It
turned out that the vanishing of the constant term in (\ref{J}) was
canonically ``forced'' by the requirement that the Monster should act
naturally on $V^{\natural}$ and on an associated algebra of vertex
operators.  This vanishing of the degree-zero subspace of
$V^{\natural}$ is actually analogous in a certain strong sense to the
absence of vectors in the Leech lattice of square-length two; the
Leech lattice is a distinguished rank-24 even unimodular (self-dual)
lattice with no vectors of square-length two.  In addition, this
vanishing of the degree-zero subspace of $V^{\natural}$ and the
absence of square-length-two elements of the Leech lattice are in turn
analogous to the absence of code-words of weight 4 in the Golay
error-correcting code, a distinguished self-dual binary linear code on
a 24-element set, with the lengths of all code-words divisible by 4.
In fact, the Golay code was used in the original construction of the
Leech lattice, and the Leech lattice was used in the construction of
$V^{\natural}$.  This was actually to be expected (if $V^{\natural}$
existed) because it was well known that the automorphism groups of
both the Golay code and the Leech lattice are (essentially) sporadic
finite simple groups; the automorphism group of the Golay code is the
Mathieu group $M_{24}$ and the automorphism group of the Leech lattice
is a double cover of the Conway group $Co_1$ mentioned above, and both
of these sporadic groups were well known to be involved in the Monster
(if it existed) in a fundamental way.  The work \cite{FLM2},
\cite{FLM5} revealed, and exploited, a new hierarchy, namely:
error-correcting codes, lattices, and vertex-operator-theoretic
structures.  The Golay code is actually unique subject to its
distinguishing properties mentioned above (proved by V. Pless
\cite{Pl}) and the Leech lattice is unique subject to its
distinguishing properties mentioned above (proved by Conway \cite{Co}
and others).  Is $V^{\natural}$ unique?  If so, unique subject to
what?  The answer to this question can be viewed as serving as a
motivation of the very notion of vertex operator algebra.  But this
uniqueness is an unsolved problem; more on this below.

After \cite{FLM2} appeared, what has been called the ``first string
theory revolution'' started in the summer of 1984, stimulated by the
work \cite{GS} of M.~Green and J.~Schwarz.  In this suddenly-active
period in string theory, the new structure $V^{\natural}$ came to be
viewed in retrospect by string theorists as an inherently {\it
string-theoretic} structure: the ``chiral algebra'' underlying the
$\mathbb{Z}_2$-orbifold conformal field theory based on the Leech
lattice.  The string-theoretic geometry is this: One takes the torus
that is the quotient of 24-dimensional Euclidean space modulo the
Leech lattice, and then one takes the quotient of this manifold by the
``negation'' involution $x \mapsto -x$, giving rise to an orbit space
called an ``orbifold''---a manifold with, in this case, a ``conical''
singularity.  Then one takes the ``conformal field theory'' (presuming
that it exists mathematically) based on this orbifold, and from this
one forms a ``string theory'' in two-dimensional space-time by
compactifying a 26-dimensional ``bosonic string'' on this
24-dimensional orbifold.  The string vibrates in a 26-dimensional
space, 24 dimensions of which are curled into this 24-dimensional
orbifold, and space-time is thus 2-dimensional in this ``toy-model''
string theory.  Such an adjunction of a two-dimensional structure is a
natural and standard procedure in string theory; 26 is the ``critical
dimension'' in bosonic string theory.  So in retrospect, the
mathematical construction \cite{FLM2} was essentially the construction
of an orbifold string theory (actually, the first example of a theory
of a string propagating on an orbifold that is not a torus).  As
discussed in the Introduction of \cite{FLM5}, some of the basic
string-theoretic papers on these aspects of orbifold string theory are
\cite{DHVW1}, \cite{DHVW2}, \cite{Ha}, \cite{DFMS}, \cite{HV},
\cite{NSV}, \cite{M}, \cite{DGH} and \cite{DVVV}.  The idea of
``orbifolding'' (as string theorists were to call it) came, in the
development of the work \cite{FLM2}, \cite{FLM5}, from the
construction of general twisted vertex operators and their algebraic
relations, including relations involving what sometimes came to be
called ``intertwining operators'' among ``twisted sectors,'' treated
in detail in \cite{FLM5} and related works discussed there.  The
construction in \cite{FLM2} also came to be viewed as a
conformal-field-theoretic structure in the sense of \cite{BPZ}, which
appeared around the same time as \cite{FLM2}.  (These ideas are all
discussed in \cite{FLM5}.)

As I mentioned at the beginning, in \cite{B1} Borcherds introduced the
axiomatic notion of {\it vertex algebra}.  This naturally extended the
relations \cite{FLM2} for the vertex operators for $V^{\natural}$ and
also other known mathematical and physical features of known vertex
operators, and these axioms turned out to be essentially equivalent to
Belavin-Polyakov-Zamolodchikov's physical axioms \cite{BPZ} for the
basic ``algebras'' of vertex operators underlying conformal field
theory.  In \cite{B1} Borcherds asserted that $V^{\natural}$ admits a
vertex algebra structure, generated by the algebraic structure
constructed in \cite{FLM2}, on which $\mathbb{M}$ (still) acts as a
symmetry group.  This assertion was proved in \cite{FLM5}, by an
(elaborate) extension of the proof of the results announced in
\cite{FLM2} (rather than by a direct use of the results announced in
\cite{FLM2}).  In \cite{FLM5}, Borcherds's definition of vertex
algebra was modified, giving the variant notion of {\it vertex
operator algebra}.  The main modification was the introduction of what
we called the ``Jacobi identity,'' discussed below.  (Actually, we had
first thought of this identity as the ``master formula'' for reasons
mentioned below, but then we decided to emphasize its analogy with the
Jacobi identity in the definition of the notion of Lie algebra.)
Another modification was the emphasis of the viewpoint that the
elements of the algebra essentially ``are'' vertex operators.  Also,
the definition of ``vertex operator algebra'' in \cite{FLM5} included
two natural grading-restriction conditions and the presence of a copy
of the Virasoro algebra, because these features naturally arose in the
construction of $V^{\natural}$.  The term ``vertex algebra'' generally
refers to any notion equivalent to Borcherds's notion in \cite{B1},
and the term ``vertex operator algebra'' generally refers to any
notion equivalent to the notion in \cite{FLM5} and also in the sequel
\cite{FHL}---that is, including the two grading restrictions---even
though, from the strictly logical point of view, the notion in
\cite{FLM5} and \cite{FHL} does not have any more ``operators'' in it
than does the notion in \cite{B1} (except in the notation).

Then in \cite{B2}, Borcherds used all this and new ideas, including
his results on generalized Kac-Moody algebras, also called Borcherds
algebras, together with certain ideas from string theory, including
the ``physical space'' of a bosonic string along with the ``no-ghost
theorem'' of R.~Brower, P.~Goddard and C.~Thorn \cite{Br}, \cite{GT},
to prove the remaining Conway-Norton conjectures for the structure
$V^{\natural}$.  What had remained to prove was that the formal series
$\sum({\rm tr}\;g|_{V^{\natural}_n})q^n$ ((\ref{gradedtraceV}) above,
but now, rather, (\ref{gradedtraceV}) for the known structure
$V^{\natural}$ instead of for a still-unknown space $V$) for the
Monster elements $g$ (or conjugacy classes) not treated in
\cite{FLM2}, \cite{FLM5}---that is, the conjugacy classes outside the
involution centralizer---were indeed the desired Hauptmoduls; the
methods of \cite{FLM2}, \cite{FLM5} did not handle these conjugacy
classes.  He accomplished this by constructing a copy of his ``Monster
Lie algebra'' from the ``physical space'' associated with
$V^{\natural}$, enlarged to a central-charge-26 vertex algebra closely
related to the 26-dimensional bosonic-string structure mentioned
above.  He transported the known action of the Monster from
$V^{\natural}$ to this copy of the Monster Lie algebra, and by using
his twisted denominator formula for this Lie algebra he proved certain
recursion formulas for the coefficients of the formal series
$\sum({\rm tr}\;g|_{V^{\natural}_n})q^n$ (that is,
(\ref{gradedtraceV}) for the known structure $V^{\natural}$), for all
Monster elements $g$.  This entailed a generalization to Borcherds
algebras of the work \cite{GL}, which had generalized B.~Kostant's
homology theorem \cite{Ko}.  The resulting recursion formulas for
$\sum({\rm tr}\;g|_{V^{\natural}_n})q^n$ agreed with the ``replication
formulas'' in \cite{CN}, satisfied by the coefficients of the
Hauptmoduls listed in \cite{CN}.  By numerically verifying that the
first few terms of some of the formal series $\sum({\rm
tr}\;g|_{V^{\natural}_n})q^n$ (for $g$ ranging through certain
elements of the involution centralizer, whose action on $V^{\natural}$
had been constructed in \cite{FLM5}) agreed with the corresponding
coefficients of the corresponding Hauptmoduls listed in \cite{CN}, he
succeeded in concluding that all the graded traces $\sum({\rm
tr}\;g|_{V^{\natural}_n})q^n$ for $V^{\natural}$ must coincide with
the formal series for the Hauptmoduls listed in \cite{CN}.

This remarkable work of Borcherds has been further illuminated in a
number of ways.  In \cite{Ju}, E. Jurisich simplified Borcherds's
argument (proving the replication formulas for the structure
$V^{\natural}$) by exploiting a certain ``large'' free Lie algebra
inside the Monster Lie algebra; this simplification is further
discussed in \cite{JLW}.  Rather than a ``Borel subalgebra'' of the
Monster Lie algebra, a certain natural ``parabolic'' subalgebra was
used, allowing the simplification.

Moreover, in \cite{CG}, C.~Cummins and T.~Gannon discovered a
conceptual proof that the replication formulas lead to the genus-zero
property.  In particular, the numerical checking in \cite{B2} using
the first few terms of the formal series constructed in \cite{FLM2},
\cite{FLM5} can essentially be bypassed.  Once one proves the
replication formulas for the action of the Monster on $V^{\natural}$
\cite{B2} (or with the shorter argument in \cite{Ju} (or \cite{JLW})),
then by \cite{CG} one knows that the ``McKay-Thompson series'' for all
the Monster elements acting on the structure $V^{\natural}$ have the
genus-zero property.  Also, the fact \cite{FLM2}, \cite{FLM5} that the
graded dimension of $V^{\natural}$ is the modular function $J(q)$,
given certain established properties of the vertex operator algebra
$V^{\natural}$, follows alternatively from a major theorem of Y.~Zhu
\cite{Z} on the modular transformation properties of the graded
dimensions of modules for suitable vertex operator algebras.

The original McKay-Thompson-Conway-Norton conjectures are conceptually
proved.  But there is also much, much more to monstrous moonshine
(some of it mentioned below).  See in particular Gannon's treatments
in \cite{Ga1} and \cite{Ga2}, which include references to many works
and surveys.

As it turned out, then, the numerology of ``monstrous moonshine'' is
much more than an astonishing relation between finite group theory and
number theory; its underlying theme is the new theory of vertex
(operator) algebras, itself the foundational structure for conformal
field theory, which is in turn the foundational structure underlying
string theory.

It's in this sense that (as I said at the beginning) vertex operator
algebra theory is inherently ``nonclassical'' in the same way in which
sporadic group theory and string theory are ``nonclassical'' in their
respective domains.

In order to motivate the precise definition of vertex (operator)
algebra, which I'll give later, I'll first repeat that there is a
vertex operator algebra (namely, $V^{\natural}$) whose symmetry group
is the Monster $\mathbb{M}$ and which implements the
McKay-Thompson-Conway-Norton conjectures relating $\mathbb{M}$ to
modular functions including $J(q)$.  But in fact, this vertex operator
algebra $V^{\natural}$ has the following three simply-stated
properties---{\it properties that have nothing at all to do with the
Monster}:

(1) {\it $V^{\natural}$, which is an irreducible module for itself
(proved in \cite{FLM5}), is its {\em only} irreducible module, up to
equivalence.}  C.~Dong \cite{D} proved this, and
C.~Dong-H.~Li-G.~Mason \cite{DLM} proved the stronger result that
every module for the vertex operator algebra $V^{\natural}$ is
completely reducible and is in particular a direct sum of copies of
itself.  {\it Thus the vertex operator algebra $V^{\natural}$ has no
more representation theory than does a field!}  (I mean a field in the
sense of mathematics, not physics.  Given a field, every one of its
modules---called vector spaces, of course---is completely reducible
and is a direct sum of copies of itself.)

(2) {\it ${\rm dim}\;V^{\natural}_0=0$}.  This corresponds to the zero
constant term of $J(q)$; while the constant term of the classical
modular function is essentially arbitrary, and is chosen to have
certain values for certain classical number-theoretic purposes, the
constant term must be chosen to be zero for the purposes of moonshine
and the moonshine module vertex operator algebra.

(3) {\it The central charge of the canonical Virasoro algebra in
$V^{\natural}$ is 24.}  ``24'' is the ``same 24'' so basic in number
theory, modular function theory, etc.  As mentioned above, this
occurrence of 24 is also natural from the point of view of string
theory.

{\it These three properties are actually ``{\it smallness}''
properties in the sense of conformal field theory and string theory.
These properties allow one to say that $V^{\natural}$ essentially
defines the smallest possible nontrivial string theory}
(cf. \cite{Ha}, \cite{Na} and the Introduction in \cite{FLM5}).
(These ``smallness'' properties essentially amount to: ``no nontrivial
representation theory,'' ``no nontrivial gauge group,'' i.e., ``no
continuous symmetry,'' and ``no nontrivial monodromy''; this last
condition actually refers to both the first and third ``smallness''
properties.)

Conversely, {\it conjecturally} \cite{FLM5}, {\it $V^{\natural}$ is
the {\em unique} vertex operator algebra with these three
``smallness'' properties (up to isomorphism).}  This conjecture turns
out to be very hard to prove (without additional strong hypotheses);
in any case, it remains unproved.  It would be the
conformal-field-theoretic analogue of the uniqueness of the Leech
lattice in sphere-packing theory and of the uniqueness of the Golay
code in error-correcting code theory, mentioned above.  Proving this
uniqueness conjecture can be thought of as the ``zeroth step'' in the
program of classification of (reasonable classes of) conformal field
theories.  M.~Tuite \cite{Tu} has related this conjecture to the
genus-zero property in the formulation of monstrous moonshine.  With
additional (strong) hypotheses assumed, uniqueness results have been
proved by Dong-Griess-C.~H.~Lam \cite{DGL} and by Lam-H.~Yamauchi
\cite{LY}.

Up to this conjecture, then, we have the following remarkable
characterization of the largest sporadic finite simple group: {\it The
Monster is the automorphism group of the smallest nontrival string
theory that nature allows, or more precisely, the automorphism group
of the vertex operator algebra with the canonical ``smallness''
properties.}  (As I mentioned above, space-time is 2-dimensional for
this ``toy-model'' string theory.  Bosonic 26-dimensional space-time
is ``compactified'' on 24 dimensions, using the orbifold construction
$V^{\natural}$; again cf. \cite{Ha}, \cite{Na} and the Introduction in
\cite{FLM5}.)  Note that (up to the conjecture) the ``smallness''
properties {\it characterize} the vertex operator algebra, but in
order to actually {\it construct} it and to construct its automorphism
group one needs the work in \cite{FLM5} or the equivalent.

This definition of the Monster in terms of ``smallness'' properties of
a vertex operator algebra provides a remarkable motivation for the
definition of the precise notion of vertex (operator) algebra.  The
discovery of string theory (as a mathematical, even if not necessarily
physical) structure sooner or later must lead naturally to the
question of whether this ``smallest'' possible nontrivial vertex
operator algebra $V^{\natural}$ exists, and the question of what its
symmetry group (which turns out to be the largest sproradic finite
simple group) is.  And on the other hand, the classification of the
the finite simple groups---{\it a mathematical problem of the
absolutely purest possible sort}---leads naturally to the question of
what natural structure the largest sporadic group is the symmetry
group of; the answer entails the development of string theory and
vertex operator algebra theory (and involves modular function theory
and monstrous moonshine as well).  The Monster, a singularly {\it
exceptional} structure---in the same spirit that the Lie algebra $E_8$
is ``exceptional,'' though $\mathbb{M}$ is {\it far} more
``exceptional'' than $E_8$---helped lead to, and helps shape, the very
{\it general} theory of vertex operator algebras.  (The exceptional
nature of structures such as $E_8$, the Golay code and the Leech
lattice in fact played crucial roles in the construction of
$V^{\natural}$, as is explained in detail in \cite{FLM3} and
\cite{FLM5}.)

Incidentally, whatever the ultimate role of string theory turns out to
be in physics, string theory is here to stay; string theory has been
``experimentally tested'' very successfully---in mathematics (whether
string theory is done by physicists or mathematicians or both), and in
many, many ways, going far beyond what I have been discussing.

The results in \cite{FLM5} include that $V^{\natural}$ is defined over
the field of real numbers, and in fact over the field of rational
numbers, in such a way that the Monster preserves the real and in fact
rational structure, and that the Monster preserves a rational-valued
positive-definite symmetric bilinear form on this rational structure.
More recent proofs that $V^{\natural}$ is a vertex operator algebra
have been found---by L.~Dolan-P.~Goddard-P.~Montague \cite{DGM}, by
Y.-Z.~Huang \cite{Hua3} and by M.~Miyamoto \cite{Mi}.  (The proof in
\cite{FLM5} is perhaps still the shortest; any complete proof must
include the full construction itself.)  Huang's proof of (the hard
part of) the vertex-operator-algebra property of $V^{\natural}$ uses
the tensor product theory for modules for a (suitable) vertex operator
algebra; I'll mention this later.  Y.~Kawahigashi and R.~Longo
\cite{KLo} have interpreted the ``orbifold'' construction of
$V^{\natural}$ in terms of algebraic quantum field theory,
specifically, in terms of local conformal nets of von Neumann algebras
on the circle.

The Monster is not the only sporadic finite simple group to which a
vertex-operator-algebraic structure has been attached.  G.~H\"ohn
\cite{Ho} has constructed a vertex-operator-superalgebraic structure
for the Baby Monster, which is involved in the Monster.  Also,
J.~Duncan \cite{Du1} has done so for the Conway group $Co_1$, and has
proved the uniqueness of the structure.  Evidence for the existence of
such a structure was given in \cite{FLM3}.  See also Borcherds-A.~Ryba
\cite{BR}.  In a remarkable development, Duncan (\cite{Du2},
\cite{Du3}) has constructed two vertex-algebraic structures for a
sporadic group {\it not} involved in the Monster, namely, the Rudvalis
group, yielding moonshine-type phenomena, including a genus-zero
property.  This supports the hope, expressed in \cite{FLM5}, that all
the sporadic groups (as well as all the other finite simple groups)
can eventually be described in vertex-algebraic terms.

So, exactly what {\it is} a vertex operator algebra?  And what {\it
are} vertex operators?  First of all, with what is now understood,
{\it vertex operators are} (or rather correspond to) {\it elements of
vertex operator algebras}, by analogy with how (for example) {\it
vectors are elements of (abstract) vector spaces}; the notion of
vector space is of course in turn defined by an axiom system.  But
before abstract vector spaces had been formalized, vectors already
``were'' something (little arrows, etc.), and this of course helped
motivate the eventual axiom system for the notion of vector space.
Here is (an oversimplified version of) what vertex operators ``already
were'':

In string theory and conformal field theory, when two (closed) strings
interact at a ``vertex,'' one has a standard picture that looks like a
``pair of pants,'' conformally equivalent to a three-punctured Riemann
sphere (after a suitable interpretation).  Such a picture is the
string-theoretic analogue of a simple Feynman diagram that looks like
the letter ``Y''---a schematic diagram for two incoming particles
interacting at a {\it ``vertex''} and producing one outgoing particle.
In the ``pair of pants,'' the {\it singularity} of the vertex in
traditional (point-particle) quantum field theory is replaced by a
{\it smooth} Riemann surface.  This allows string theory to avoid the
``ultraviolet divergences'' in point-particle quantum field theory.
In string theory and conformal field theory, such ``Riemann-surface
vertex diagrams'' get ``represented'' by ``vertex operators'' acting
on suitable infinite-dimensional vector spaces; vertex operators
``describe'' the particle (or rather, string) interactions in a given
conformal field theory model.  Geometric relations among
Riemann-surface ``diagrams'' are reflected by algebraic and analytic
relations among vertex operators.

I'll next give a concrete example of a vertex operator, as it arose in
mathematics:

For certain mathematical reasons, with R.~Wilson \cite{LW1} we focused
on the ``philosophical'' problem of trying to construct the affine
Kac-Moody Lie algebra
\begin{equation}\label{sl2hat}
A_1^{(1)} = \widehat{\mathfrak{s}\mathfrak{l}(2)}=\mathfrak{s}\mathfrak{
l}(2)\otimes \mathbb{C}[t,t^{-1}]\oplus\mathbb{C}c
\end{equation}
(the ``affinization'' of the Lie algebra $\mathfrak{s}\mathfrak{l}(2)$ of
2-by-2 matrices of trace 0), with Lie brackets given by
\[
[a\otimes t^m,b\otimes t^n]=[a,b]\otimes t^{m+n}+{\rm
tr}(ab)m\delta_{m+n,0}c
\]
for $a,b \in \mathfrak{s}\mathfrak{l}(2)$ and $m,n \in \mathbb{Z}$ and
\[
[c,\widehat{\mathfrak{s}\mathfrak{l}(2)}]=0,
\]
as some kind of ``concrete'' Lie algebra of (as-yet unknown) operators
on some kind of ``concrete'' space.

In fact, our main reasons for formulating and trying to solve this
problem were that we already knew some pieces of evidence, stemming
{}from joint work with S.~Milne \cite{LM} and with A.~Feingold
\cite{FL}, that such a construction might eventually shed light on the
classical Rogers-Ramanujan combinatorial identities.  One of these two
identities states that the number of partitions of a nonnegative
integer $n$ into parts congruent to 1 or 4 mod 5 equals the number of
partitions of $n$ into parts whose successive differences are at least
2, and the other of these identities states that the number of
partitions of a nonnegative integer $n$ into parts congruent to 2 or 3
mod 5 equals the number of partitions of $n$ into parts whose
successive differences are at least 2 and such that the smallest part
is at least 2.  These two theorems have a long and interesting history
and are highly nontrivial; cf. \cite{A1}.  When these two identities
are written in their original, classical, generating-function form
(cf. \cite{A1}), each of them asserts the equality of two formal power
series ($q$-series)---one of them a formal infinite product in $q$ and
the other a formal infinite sum.  The work \cite{LM}, which used the
Weyl-Kac character formula \cite{Ka}, showed that the product sides of
the two Rogers-Ramanujan identities had something interesting to do
with standard ($=$ integrable highest weight) $\widehat{\mathfrak{
s}\mathfrak{l}(2)}$-modules of levels 1 and 3; the ``level'' is the
scalar by which the central element $c$ in (\ref{sl2hat}) acts.  It
seemed natural to try to ``construct'' these standard modules somehow,
starting with the level 1 standard modules (the ``basic'' modules).
The hope was to try to ``discover'' the sum sides of the
Rogers-Ramanujan identities, somehow, in the level 3 standard modules.
The Rogers-Ramanujan identities had been proved many times, but the
question now was: What do the sum sides of the Rogers-Ramanujan
identities ``count,'' in this new context?  (Classically, they count
partitions satisfying the difference-two condition.)

Here is the result of the work \cite{LW1} (expressed in different, but
equivalent, notation): Consider the (commutative associative) algebra
\[
S=\mathbb{C}[y_{\frac{1}{2}},y_{\frac{3}{2}},y_{\frac{5}{2}},\dots]
\]
of polynomials in the formal variables $y_n$,
$n=\frac{1}{2},\frac{3}{2},\dots$.  Form the expression
\[
Y(x)={\rm exp}\left(\sum_{n=\frac{1}{2},\frac{3}{2},\frac{5}{2},\dots}
\frac{y_n}{n}x^n\right) {\rm
exp}\left(-2\sum_{n=\frac{1}{2},\frac{3}{2},\frac{5}{2},\dots}
\frac{\partial}{\partial y_n}x^{-n}\right),
\]
 where ``{\rm exp}'' is the formal exponential series and $x$ is
another formal variable commuting with the $y_n$'s.  The $y_n$'s
(understood as multiplication operators on $S$) can be thought of as
``creation operators'' and the $\frac{\partial}{\partial y_n}$'s as
``annihilation operators,'' acting on the ``Fock space'' $S$, using
some terminology from quantum field theory.  Together with the
identity operator on the space $S$, they span an
(infinite-dimensional) Heisenberg Lie algebra acting on $S$; the
commutators among these operators are the classical Heisenberg
commutation relations, on infinitely many generators.  The operator
$Y(x)$ is a well-defined formal differential operator in infinitely many
formal variables, including the extra variable $x$.  Viewing $Y(x)$ as a
generating function with respect to the formal variable $x$, we write
\[
Y(x)=\sum_{j \in {\frac{1}{2}}\mathbb{Z}} A_j x^{-j},
\]
thus giving a family of (well-defined) linear operators $A_j$, $j \in
{\frac{1}{2}}\mathbb{Z}$, acting on $S$.  Each $A_j$ can be computed, as
a certain formal differential operator, in the form of an infinite sum
of products of multiplication operators with partial differentiation
operators, multiplied in this order; this infinite sum actually acts
as a finite sum when applied to any given element of the space $S$.
The explicit expression for each $A_j$ is in fact complicated, and
while one {\it can} write it down explicitly, one does not want to
have to do this, although in our original work we did in fact find
these explicit formal differential operators $A_j$ ``directly''; it
was only after the fact that we realized that if we added up {\it all}
of these complicated operators $A_j$ and thus formed their generating
function as above, then all of these operators $A_j$ could be
described by the single product of exponentials $Y(x)$, which looked
much simpler than any of the individual operators $A_j$.

The main point of this was:

\begin{theorem}\label{theoremsl2} \cite{LW1}  The operators
\[
1, \;\;y_n,\;\;\frac{\partial}{\partial y_n}\;{\rm \Big(}
n=\frac{1}{2},\frac{3}{2},\frac{5}{2},\dots{\rm \Big)}\;\; {\rm
and}\;\;A_j\;{\rm \Big(} j \in \frac{1}{2}\mathbb{Z}{\rm \Big)}
\]
(1 is the identity operator on the space $S$) span a Lie algebra of
operators acting on $S$, that is, the commutator of any two of these
operators is a (finite) linear combination of these operators, and
this Lie algebra is a copy of the affine Lie algebra $\widehat{\mathfrak{
s}\mathfrak{l}(2)}$.
\end{theorem}

This operator $Y(x)$ turned out to be a variant of the vertex
operators that had arisen in string theory, as H.~Garland pointed out
(although that is not how we had found it).  It turned out that vertex
operators and symmetry are closely related.  Instances of this had
already been discovered in physics, including the works \cite{H} and
\cite{BHN}.

With what is now known, this operator $Y(x)$ is an example of a
``twisted vertex operator'' (a vertex operator appearing in a twisted
module for a vertex operator algebra).  This particular vertex
operator construction happens to be among the (many) ingredients
playing crucial roles in the construction of $V^{\natural}$, as are
the untwisted vertex operator constructions of I.~Frenkel-V.~Kac
\cite{FK} and G.~Segal \cite{S1} as well.  These and other vertex
operator constructions also enter into a variety of other, very
different, mathematical problems.  For instance, the operator $Y(x)$
above was interpreted by E.~Date, M.~Kashiwara and T.~Miwa \cite{DKM}
to be precisely the infinitesimal B\"acklund transformation for the
Korteweg-de Vries hierarchy of differential equations in soliton
theory; in Hirota's bilinear formalism, $Y(x)$ generates the
multi-soliton solutions.

In Theorem \ref{theoremsl2}, the $\widehat{\mathfrak{s}\mathfrak{
l}(2)}$-module that is constructed using the twisted vertex operator
$Y(x)$ is one of the two basic (level 1 standard) $\widehat{\mathfrak{
s}\mathfrak{l}(2)}$-modules.  The other one of the two basic modules is
constructed in exactly the same way, but with each operator $A_j$ in
Theorem \ref{theoremsl2} replaced by $-A_j$ (or equivalently, with the
generating-function $Y(x)$ of the operators $A_j$ replaced by its
negative $-Y(x)$).

Further work with Wilson led to the construction of structures we
called ``$Z$-algebras'' \cite{LW2}--\cite{LW4}, which provided a
vertex-operator-theoretic interpretation and proof the
Rogers-Ramanujan identities (mentioned above), in the following way
(very briefly):

The (higher-level) standard $\widehat{\mathfrak{s}\mathfrak{l}(2)}$-modules
of level $k>1$ can be constructed inside the tensor products of $k$
copies of basic modules, but it was an open problem to construct these
higher-level modules ``concretely,'' by exhibiting natural bases of
them.  The work \cite{LM} led to a natural conjecture that the
Rogers-Ramanujan identities should ``take place'' inside the level 3
standard $\widehat{\mathfrak{s}\mathfrak{l}(2)}$-modules, with a structure
that looks like a copy of the basic $\widehat{\mathfrak{s}\mathfrak{
l}(2)}$-module somehow ``factored out,'' as a tensor factor.  It
turned out that the Heisenberg Lie subalgebra of $\widehat{\mathfrak{
s}\mathfrak{l}(2)}$ entering into Theorem \ref{theoremsl2} acts
completely reducibly on each standard module $L$, and that the
``vacuum space'' $\Omega_L$ in $L$ for the action of this Heisenberg
subalgebra (the subspace of $L$ annihilated by all the ``annihilation
operators'' mentioned above) has an easily-computed graded dimension.
This vacuum space $\Omega_L$ implemented the desired ``factoring
out,'' and in case the level of $L$ is 3, the $q$-series that are the
graded dimensions of the spaces $\Omega_L$ are exactly the product
sides of the two Rogers-Ramanujan identities.  The next problem
appeared to be to find a basis of the vacuum space $\Omega_L$ (for a
level 3 standard module) that would exhibit the graded dimension of
$\Omega_L$ as the {\it sum} side of the corresponding Rogers-Ramanujan
identity---the classical $q$-series generating function of the number
of partitions of $n$ satisfying one of the two difference-two
conditions mentioned above.

Using the type of structure involved in Theorem \ref{theoremsl2}, in
\cite{LW2}--\cite{LW4} we eventually constructed certain operators
that commute with the action of the Heisenberg Lie subalgebra of
$\widehat{\mathfrak{s}\mathfrak{l}(2)}$, acting on suitable modules, and we
called these operators ``$Z$-operators'' (``$Z$'' referring to the
centralizing of this subalgebra).  A typical $Z$-operator is a
generating function of the shape
\[
Z(x)=\sum_{j \in \mathbb{Z}} Z_j x^{-j},
\]
where each $Z_j$ is an operator of ``degree'' $j$ that (because of the
centralizing property) preserves $\Omega_L$.  The goal seemed to be to
prove that the vacuum space $\Omega_L$, for a level 3 standard module,
has a basis of the form
\begin{equation}\label{monomials}
Z_{j_1}Z_{j_2} \cdots Z_{j_n}, \;\;\; j_i < 0,
\end{equation}
with
\begin{equation}\label{differencetwo}
j_1 \leq j_2 -2, \;\;  j_2 \leq j_3 -2, \dots, j_{n-1} \leq j_n -2,
\end{equation}
and with $j_n \leq -2$ as well, in the case of the second of the
Rogers-Ramanujan identities.  (The subscripts $j_i$ play the role of
the negatives of the parts in a partition of a nonnegative integer.)
This difference-two condition on the indices of such a basis monomial
would exhibit the graded dimension of $\Omega_L$ as the sum side of a
Rogers-Ramanujan identity, and this would solve the problem.  But
achieving this would require the ``straightening'' of monomials
(\ref{monomials}) to obtain the inequalities (\ref{differencetwo}) on
the indices.  Such ``straightening'' could be accomplished if there
were good enough algebraic relations involving the generating function
$Z(x)$, such as perhaps a commutator formula for $[Z(x_1),Z(x_2)]$.

But there {\it is} no such commutator formula.  Instead, what turned
out to be possible was the construction of ``generalized commutation
relations'' of the shape
\begin{equation}\label{generalizedcommrels}
A(x_1,x_2)Z(x_1)Z(x_2) - B(x_1,x_2)Z(x_2)Z(x_1) = C(x_1,x_2),
\end{equation}
where $A(x_1,x_2)$ and $B(x_2,x_1)$ are suitable formal expansions of
suitable formal algebraic functions, and $C(x_1,x_2)$ is some operator
that is ``simpler than'' both $Z(x_1)Z(x_2)$ and $Z(x_2)Z(x_1)$.  The
term ``generalized'' refers to the presence of the formal algebraic
functions $A$ and $B$; in the case when these are 1, then one of
course has ordinary commutation relations.  Also, there are
``generalized anticommutation relations,'' of a generally
still-more-complicated shape, that involve arbitrary numbers of
generating functions $Z(x_i)$ in general.  {\it Each} of these
generalized commutation and anticommutation relations can be viewed as
the generating function of an {\it infinite family of relations} among
monomials in the operators $Z_j$; each such relation among such
monomials is of the following form: A (well-defined) formal infinite
linear combination of monomials in the $Z_j$'s, with coefficients
coming from the coefficients of the formal algebraic functions such as
$A$ and $B$, is equated with a ``simpler'' expression involving the
$Z_j$'s.  (Again, if the formal algebraic functions $A$ and $B$ are 1
in (\ref{generalizedcommrels}), then each such relation among the
$Z_j$'s is a commutation relation of the following form:
$[Z_{j_1},Z_{j_2}]$ equals a simpler expression such as perhaps a
multiple of a single $Z_j$ or a scalar.)  These various types of
generalized commutation and anticommutation relations entering into
the solution of the present problem are detailed in
\cite{LW2}--\cite{LW4}.  What they were used for was to ``straighten''
monomials (\ref{monomials}) to obtain the difference-two condition
(\ref{differencetwo}) on the subscripts, in the case of the the level
3 standard $\widehat{\mathfrak{s}\mathfrak{l}(2)}$-modules.  This showed
that the monomials (\ref{monomials}) satisfying the difference-two
condition {\it span} the vacuum space $\Omega_L$.  Much more work was
needed to prove that these difference-two monomials are also {\it
linearly independent}, thus proving that they form a basis.  Once
these difference-two monomials were proved to form a basis, this
vertex-operator interpretation and proof of the Rogers-Ramanujan
identities was complete.

I had mentioned generalized commutation relations when I was
discussing how certain such relations among certain vertex operators
gave a natural ``commutative affinization'' of the (enlargement of)
the Griess algebra acting on $V^{\natural}$, in \cite{FLM2},
\cite{FLM5}.  In that situation, the formal algebraic functions $A$
and $B$ in (\ref{generalizedcommrels}) are particularly (but
deceptively!)  simple:
\begin{equation}
A(x_1,x_2) = B(x_2,x_1) = x_1 - x_2.
\end{equation}
That is, the relations yielding the commutative affinization of the
196884-dimensional enlargement ${\mathcal B}$ of the Griess algebra in
$V^{\natural}$ are of the shape
\begin{equation}\label{crossbracket}
(x_1-x_2)[Z_1(x_1),Z_2(x_2)] = {\rm simpler},
\end{equation}
where in this formula, the generating functions $Z_i(x_j)$ now refer to
the relevant vertex operators in \cite{FLM2}, \cite{FLM5}.  These
relations, which we called ``cross-bracket'' relations in \cite{FLM2},
\cite{FLM5}, exhibit the desired commutative affinization, once one
has constructed the space $V^{\natural}$.

The process discussed above for the level 3 standard $\widehat{\mathfrak{
s}\mathfrak{l}(2)}$-modules was extended to all the standard
$\widehat{\mathfrak{s}\mathfrak{l}(2)}$-modules in \cite{LW2}--\cite{LW4},
and the result was a corresponding vertex-algebraic interpretation of
a known family of generalizations of the Rogers-Ramanujan identities,
which had been discovered by B.~Gordon, G.~Andrews and D.~Bressoud.
Many of these identities are treated in \cite{A1}.  The case of the
level 2 standard $\widehat{\mathfrak{s}\mathfrak{l}(2)}$-modules actually
gave certain infinite-dimensional Clifford algebras, and
correspondingly, a ``difference-one'' condition, familiar in a natural
basis of an exterior algebra, rather than a ``difference-two
condition.''  It was for this reason that we thought of the general
phenomenon, arising now for all the standard $\widehat{\mathfrak{s}\mathfrak{
l}(2)}$-modules, as the emergence of a new ``generalized Pauli
exclusion principle'' for a natural family of operators generalizing
classical Clifford-algebra operators that are familiar in quantum
mechanics and quantum field theory and producing fermionic particles;
starting at level 3, fermionic (``difference-one'') statistics changed
into ``difference-two statistics,'' and for the levels greater than 3,
``difference-two-at-a-distance statistics,'' reflecting the sum sides
of the Gordon-Andrews-Bressoud identities.  (See
\cite{LW2}--\cite{LW4}.)

However, we were unable to prove the linear independence of the
relevant monomials, analogous to those (for the level-3 case) in
(\ref{monomials}), (\ref{differencetwo}), for the levels greater than
3.  This problem was solved by A.~Meurman and M.~Primc \cite{MP1}, who
thus provided a vertex-algebraic proof of the Gordon-Andrews-Bressoud
identities beyond the case of the Rogers-Ramanujan identities.  Also,
C.~Husu \cite{Hus1} discovered an elegant ``multi-Jacobi-identity''
interpretation and proof of the complicated generalized
anticommutation relations in \cite{LW2}, \cite{LW4}.

The $Z$-algebra viewpoint in \cite{LW2}--\cite{LW4} was used by
K.~Misra (\cite{Mis1}--\cite{Mis4}), M.~Mandia \cite{Ma}, C.~Xie
\cite{X}, S.~Capparelli \cite{Ca1} and M.~Bos-K.~Misra \cite{BosMis}
to construct difference-two-type bases for the vacuum spaces
$\Omega_L$ for certain standard modules for a range of affine Lie
algebras, giving still more interpretations and proofs of the
Rogers-Ramanujan and Gordon-Andrews-Bressoud identities, and Bos
\cite{Bos} proved that in fact the complete list of such occurrences
of the Rogers-Ramanujan identities consists of certain low-level
standard modules for the affine Kac-Moody algebras $A_1^{(1)}$,
$A_2^{(1)}$, $A_2^{(2)}$, $A_7^{(2)}$, $C_3^{(1)}$, $F_4^{(1)}$ and
$G_2^{(1)}$; all these cases are covered by the papers mentioned.

One goal was to discover new identities with these ideas.  In his
Ph.D. thesis research at Rutgers, S.~Capparelli set out to construct
$Z$-algebra bases of the vacuum spaces $\Omega_L$ for the standard
$A_2^{(2)}$-modules of level 3, and he succeeded in obtaining a
construction, and a pair of partition identities, but only as a
conjecture, because he had not yet completed his $Z$-algebra proof of
the linear independence of his spanning sets of $Z$-operator monomials
in the vacuum spaces.  (He completed this in \cite{Ca2}, and
M.~Tamba-C.~Xie \cite{TX} did so as well.)  Meanwhile, in a talk at
the Rademacher Centenary Conference in 1992, I mentioned Capparelli's
still-conjectured identities, or rather, one of them, in a survey of
the research program.  After confirming that the identity was indeed
new, Andrews proved it, before the conference had ended (see
\cite{A2}), and soon afterward, K.~Alladi, Andrews and Gordon
formulated and proved a refinement and generalizations of it
\cite{AAG}.

A central goal of these ideas was to find interesting new structure.
The $Z$-algebra constructions mentioned so far are based on a certain
twisting, the same twisting as in \cite{LW1}.  In joint work with
Primc \cite{LP}, we developed $Z$-algebras in the untwisted case, and
using these and related ideas in this analogous but still surprisingly
different setting, we constructed combinatorial bases exhibiting
``difference-two-at-a-distance generalized fermionic statistics,'' and
corresponding (``fermionic'') character formulas, for the higher-level
standard $\widehat{\mathfrak{s}\mathfrak{l}(2)}$-modules.  In a sequel
\cite{MP3} to \cite{MP1}, Meurman and Primc discovered and exploited
new structure related to \cite{LP}.

A.~Zamolodchikov-and V.~Fateev introduced ``nonlocal parafermion
currents'' in \cite{ZF1} and a twisted analogue in \cite{ZF2}, and
these conformal-field-theoretic constructions turned out to be
essentially equivalent to the untwisted and twisted $Z$-algebra
constructions in \cite{LP} and \cite{LW2}, respectively; see
\cite{DL1}, \cite{DL2} for the untwisted case, and \cite{Hus1} for the
twisted case.  The ``nonlocality'' in Zamolodchikov-Fateev's
terminology refers to the fact that in the notation above, the formal
algebraic functions $A(x_1,x_2)$ and $B(x_1,x_2)$ in
(\ref{generalizedcommrels}) are (formally) multiple-valued, and the
term ``parafermion'' refers to the generalization of fermionic
statistics mentioned above, or more precisely, to the reflection of
such statistics in the form of the formal algebraic functions $A$ and
$B$.  In the terminology of \cite{DL2}, the $Z$-algebra structures
generate certain examples of ``generalized vertex algebras'' and
``abelian intertwining algebras,'' whose main axiom is a ``generalized
Jacobi identity,'' incorporating the relevant formally-multiple-valued
formal algebraic functions.  Notions essentially equivalent to
generalized vertex algebras were also introduced by
A.~Feingold-I.~Frenkel-J.~Ries \cite{FFR} and G.~Mossberg \cite{Mos}.
The introductory material in \cite{DL2} includes discussions of these
developments.

One of the fundamental methods classically used to study and prove the
Rogers-Ramanujan identities was the Rogers-Ramanujan recursion, as
explained in \cite{A1}.  But this recursion did not arise in any of
the vertex-operator work on partition identities that I have
mentioned.  More recently, in joint work work with Capparelli and
A.~Milas (\cite{CapLM1}, \cite{CapLM2}), using vertex operator algebra
theory in the context of ``principal subspaces'' (\cite{FS1},
\cite{FS2}) of standard modules, we have been able to incorporate
these into the theory, as well as the more general Rogers-Selberg
recursions, satisfied by the sum sides of the Gordon-Andrews
identities.  This leads to new questions, being addressed in work of
C.~Calinescu \cite{Cal1}, \cite{Cal2} and joint work with Calinescu
and Milas \cite{CalLM}.  The main theme here is to use intertwining
operators (which I'll mention below) among modules for vertex
(operator) algebras to find new structure.  In fact, perhaps the main
theme of all the work I've mentioned on partition identities is to use
known and sometimes new identities as clues to look for interesting
new structure.

Now I'll give the definition of the notion of vertex operator algebra.
The considerations I've mentioned (and related additional ones) led to
the following variant (\cite{FLM5} and \cite{FHL}) of Borcherds's
definition of the notion of vertex algebra; in this definition we use
commuting independent formal variables $x$, $x_1$, $x_2$, etc.:

\begin{definition}\label{vertexoperatoralgebra}
{\rm A {\it vertex operator algebra} consists of a $\mathbb{Z}$-graded vector
space
\[
V=\bigoplus_{n\in \mathbb{Z}}V_{(n)}
\]
(this grading, by {\it conformal weights}, is ``shifted'' from the
grading we have already been using on $V^{\natural}$; that grading is
adapted to the modularity properties of the generating functions,
including the $J$-function, that we have mentioned, while the present
grading is adapted to the action of the Virasoro algebra on the vertex
operator algebra, as we mention below) such that
\[
\mbox{\rm dim }V_{(n)}<\infty\;\;\mbox{\rm for}\; n \in \mathbb{Z},
\]
\[
V_{(n)}=0\;\;\mbox{\rm for} \;n\; \mbox{\rm sufficiently negative},
\]
equipped with a linear map (the {\it ``state-field correspondence''})
\begin{eqnarray*}
V&\to&(\mbox{\rm End}\; V)[[x, x^{-1}]]\nonumber \\
v&\mapsto& Y(v, x)={\displaystyle \sum_{n\in\mathbb{Z}}}v_{n}x^{-n-1}
\end{eqnarray*}
(where $v_{n}\in {\rm End}\;V$, the algebra of linear operators on
$V$), $Y(v, x)$ denoting the {\it vertex operator associated with $v$}
(the letter ``$Y$'' happens to look like the vertex Feynman diagram
that we mentioned above), and equipped also with two distinguished
homogeneous vectors ${\bf 1}\in V_{(0)}$ (the {\it vacuum vector}) and
$\omega \in V_{(2)}$ (the {\it conformal vector}).  The axioms are:
For $u, v \in V$,
\[
u_{n}v=0 \;\; {\mbox {\rm for}} \;n\; {\mbox {\rm sufficiently large}}
\]
(the {\it truncation condition});
\[
Y({\bf 1}, x) = {\mbox {\rm the identity operator on}} \; V;
\]
the {\it creation property}:
\[
Y(v, x){\bf 1} \; {\mbox {\rm has no pole in}} \; x \; {\mbox {\rm and
its constant term is the vector}} \; v
\]
(which implies that the state-field correspondence is one-to-one); the
{\it Virasoro algebra relations}:
\[
[L(m), L(n)]=(m-n)L(m+n)+\frac{m^{3}-m}{12}\delta_{n+m,0}c
\]
for $m, n \in \mathbb{Z}$, where $c\in \mathbb{C}$ (the {\it central
charge}) and where
\[
Y(\omega, x)=\sum_{n\in\mathbb{Z}}L(n)x^{-n-2}
\]
(i.e., the conformal vector generates the Virasoro algebra);
\[
L(0)v=nv
\]
for $n \in \mathbb{Z}$ and $v \in V_{(n)}$;
\[
{\displaystyle\frac{d}{dx}}Y(v, x)=Y(L(-1)v, x);
\]
and finally, the main axiom (the vast bulk of the definition of vertex
operator algebra), the {\it Jacobi identity}: Writing
\begin{equation}\label{delta}
\delta(x)=\sum_{n\in \mathbb{Z}}x^n,
\end{equation}
the ``formal delta function'' (this really {\it is} a natural analogue
of the Dirac delta function), we have:
\begin{eqnarray}\label{Jacobi}
\lefteqn {\displaystyle x^{-1}_0\delta \left( \frac{x_1-x_2}{x_0}\right)
Y(u,x_1)Y(v,x_2)}\nonumber\\
&&\hspace{2em}- x^{-1}_0\delta \left( \frac{-x_2+x_1}{x_0}\right)
Y(v,x_2)Y(u,x_1)\nonumber \\
&&{\displaystyle = x^{-1}_2\delta \left( \frac{x_1-x_0}{x_2}\right)
Y(Y(u,x_0)v,x_2)},
\end{eqnarray}
where each binomial, such as for example $(x_1-x_2)^{-n}$, occurring
in the expansions here is understood to be expanded in nonnegative
integral powers of the {\it second} variable; the truncation condition
insures that all the expressions here are well defined.  The identity
(\ref{Jacobi}) is actually the generating function of an infinite list
of identities among operators on $V$, obtained by equating the
coefficients of all the monomials in $x_0$, $x_1$ and $x_2$.}
\end{definition}

The ``physical'' terms mentioned in this definition reflect the
well-understood relation between this algebraic notion and the notion
of ``chiral algebra'' in conformal field theory.  The ``language'' of
this definition, in particular, of the Jacobi identity, is ``formal
calculus,'' which is discussed in detail in \cite{FLM5}, \cite{FHL}
and \cite{LL}, for instance.  We first learned about the viewpoint of
the formal delta function (\ref{delta}) from Garland; see also
\cite{DKM}.

Formal calculus is completely different from the algebra of formal
power series.  Formal power series (with the powers of the formal
variable(s) all nonnegative), or formal Laurent series with the powers
of the formal variable(s) truncated from below, form rings.  They can
be multiplied.  But the formal series in the ``formal calculus'' of
vertex operator algebra theory are ``doubly infinite''; the powers of
the formal variables in a given formal series can be arbitrary
negative {\it or} positive integers, as is the case in (\ref{delta})
and (\ref{Jacobi}), and sometimes even more generally, the powers of
the formal variables in a given formal series must be allowed to be
arbitrary rational or even {\it complex} numbers.  These general kinds
of ``doubly infinite formal Laurent series'' and so on certainly do
{\it not} form rings.  Yet there is a highly-developed algebra
(``formal calculus'') for such formal series, and indeed, one needs
this in vertex (operator) algebra theory and in particular in the
representation theory (for instance, in the tensor product theory for
module categories for a suitable vertex operator algebra; I'll mention
this later).  A typical pattern in the theory is: First one does
formal calculus, and then one (suitably and systematically)
specializes the formal variables to complex variables and proves
(analytic) convergence, for those results that need to be formulated
analytically.

Now that we have the precise definition of the notion of vertex
operator algebra, we have a canonical definition of the Monster
without reference to finite group theory: It is the symmetry group of
the (conjecturally) unique vertex operator algebra (a structure
satisfying the Jacobi identity and the ``relatively minor'' axioms)
having the three ``smallness'' properties mentioned above.  Thus the
Jacobi identity plus a few words determine the largest sporadic finite
simple group.

Notice that there is a certain resemblance between the ``shape'' of
the Jacobi identity (\ref{Jacobi}) and the (schematic) generalized
commutation relation (\ref{generalizedcommrels}); the three-variable
formal-delta-function expressions in the left-hand side of
(\ref{Jacobi}) are analogous to the formal expressions $A(x_1,x_2)$
and $B(x_1,x_2)$ in (\ref{generalizedcommrels}).  This is actually
more than a coincidence.  In \cite{B1}, Borcherds had discovered
identities equivalent, in retrospect, to certain consequences of
(\ref{Jacobi}).  These included in particular a formula for the
commutator $[Y(u,x_1),Y(v,x_2)]$, which can be obtained from
(\ref{Jacobi}) by taking the (formal) residue with respect to the
formal variable $x_0$ (that is, the coefficient of $x_0^{-1}$) of each
of the two sides of (\ref{Jacobi}) and equating the results, as well
as a formula for $Y(Y(u,x_0)v,x_2)$, which can be obtained by taking
the residue with respect to the formal variable $x_1$ of each of the
two sides of (\ref{Jacobi}) and equating the results.  In \cite{B1},
these and other identities for a vertex algebra were actually
expressed in component form (that is, in terms of the component
operators $v_n$ of the vertex operators $Y(v,x)$) instead of
generating-function form, but they could be recast in
generating-function form.  Also, Borcherds used the second one of
these two identities, together with certain other axioms, not
including his commutator formula, is his definition of ``vertex
algebra.''  When expressed, in retrospect, using the formal variables
of (\ref{Jacobi}), each of these two identities involves only two of
the three formal variables $x_0$, $x_1$, $x_2$ (the third variable
being eliminated by the residue procedure).

Now, the identity of the shape (\ref{crossbracket}) that was needed
for the commutative affinization, using ``cross-brackets,'' in
\cite{FLM2}, \cite{FLM5} suggested the formulation of a {\it family}
of generalized commutation relations, with the left-hand sides of the
form
\begin{equation}\label{nthcrossbracket}
(x_1-x_2)^n Z_1(x_1)Z_2(x_2) - (-x_2+x_1)^n Z_2(x_2)Z_1(x_1)
\end{equation}
for all $n \in \mathbb{Z}$ (where the binomial expansion convention
mentioned just after (\ref{Jacobi}) is being used here), so that the
case $n=0$ would give a commutator formula and the case $n=1$ would
give an expression for the cross-bracket of type (\ref{crossbracket}).
In \cite{FLM5} we were in fact able to do this, and the most natural
way to formulate the resulting infinite family of identities was to
put all these identities into a single generating function, using a
new formal variable $x_0$; that is, there was now a formula for
\begin{equation}\label{allcrossbrackets}
\sum_{n \in \mathbb{Z}}x_0^{-n-1}(x_1-x_2)^n Z_1(x_1)Z_2(x_2)
-\sum_{n \in \mathbb{Z}}x_0^{-n-1}(-x_2-x_1)^n Z_2(x_2)Z_1(x_1)
\end{equation}
(where the operators called $Z_i$ here are now the vertex operators
entering into the construction of $V^{\natural}$).  {\it But
(\ref{allcrossbrackets}) is exactly the left-hand side of
(\ref{Jacobi})} (with the operators $Z_i$ now playing the roles of the
vertex operators in (\ref{Jacobi})).  In other words, (\ref{Jacobi})
is the generating function of an infinite family of generalized
commutation relations.  Also, we were able to prove in \cite{FLM5}
that this generating function actually equaled the right-hand side of
(\ref{Jacobi}), for all the vertex operators involved in
$V^{\natural}$.  Thus for us, in \cite{FLM5} the idea for formulating,
and proving, the Jacobi identity arose from the axioms and properties
that Borcherds wrote in \cite{B1} together with the
generalized-commutation-relation idea.  Notice that the three
expressions in (\ref{Jacobi}) are analogous, and even quite similar,
to one another.  In fact in \cite{FHL}, an explicit three-variable
symmetry of the Jacobi identity was formulated and analyzed.  If one
takes the residue with respect to any of the three formal variables in
(\ref{Jacobi}), the three-variable symmetry is ``broken,'' and
correspondingly, with (\ref{Jacobi}) not yet known to us, it seemed
natural to try to construct the identity (\ref{Jacobi}), completing
this natural symmetry.  The identity (\ref{Jacobi}) exhibits ``all''
the information.  It was for these reasons that we were initially
thinking of (\ref{Jacobi}) as the ``master formula.''  (As we said in
the Introduction of \cite{FLM5}, Borcherds informed us that he too
found this identity, and in fact it is implicit in \cite{B1}.  It is
also implicit in the work of physicists, as we have been discussing.)

But we decided instead in \cite{FLM5} to call (\ref{Jacobi}) the
Jacobi identity because it is analogous to the Jacobi identity in the
definition of Lie algebra: For $u$ and $v$ in a Lie algebra,
\[
({\rm ad}\; u)({\rm ad}\; v)-({\rm ad}\; v)({\rm ad}\; u)={\rm
ad}(({\rm ad}\; u)v),
\]
where ${\rm ad}\; u$ is the operation of left-bracketing with the
element $u$.  While the variables in (\ref{Jacobi}) must be understood
as formal and {\it not} complex variables (for instance, note that in
the formal delta function (\ref{delta}) itself, the formal variable
cannot be specialized to a complex variable; the resulting
doubly-infinite formal series converges nowhere), it is in fact
possible (and very important) to specialize the formal variables to
complex variables in certain systematic, and subtle, ways.  This
process is crucial, in particular, for making the connection between
the notion of vertex operator algebra and the notion of chiral algebra
in conformal field theory, and it is also crucial for mathematical
reasons.  In fact, because of the fundamental complex-analytic
geometry implicit in (\ref{Jacobi}), suitably interpreted, we also
called (\ref{Jacobi}) the ``Jacobi-Cauchy identity'' in \cite{FLM5},
for reasons explained in detail in the Appendix of \cite{FLM5}.

The notion of vertex operator algebra, then, is indeed deeply
analogous to the notion of Lie algebra, and {\it is actually the
``one-complex-dimensional analogue'' of the notion of Lie algebra}
(which is the corresponding ``one-real-dimensional'' notion, in this
sense); this statement can be made precise using the language and
viewpoint of operads, which I will mention below in connection with
Huang's work incorporating the geometry underlying the Virasoro
algebra into the structure.

In addition to being the generating function of an infinite family of
generalized-commutation-relation identities, (\ref{Jacobi}), as we
mentioned above, is the generating function of an infinite family of
(generally highly-nontrivial) identities for the component operators
$v_n$ of the vertex operators in a vertex operator algebra $V$, one
identity for each monomial in the three formal variables $x_0$, $x_1$
and $x_2$.  Each identity in this infinite list involves binomial
coefficients, coming from the three formal delta-function expressions.
It is the {\it generating function form of these identities (namely,
the Jacobi identity)} that is the natural analogue of the Jacobi
identity in the definition of Lie algebra.  Even more basically,
$Y(u,x)$ is itself the generating function of the infinite family of
operators $u_n$ acting on $V$, as we saw in the concrete example of
the (twisted) vertex operator $Y(x)$ above (acting on the space $S$,
in that case; as we emphasized, the generating function is much easier
to work with than the individual operators $A_j$).  Moreover, the
single ``$x$-parametrized product operation'' taking the ordered pair
$(u,v)$, $u,v \in V$, to the generating function $Y(u,x)v$ can
certainly be thought of as specifying an infinite family of
nonassociative product operations $u_n v$ on $V$, for $n \in 
\mathbb{Z}$, corresponding to the powers of $x$.  In very special cases, the
``component identities'' of the Jacobi identity include the relations
defining affine Lie algebras; the Virasoro algebra; the
infinite-dimensional ``affinization'' of the (modified) Griess algebra
${\mathcal B}$; and a vast array of other remarkable algebraic
structures.  In ``formal calculus,'' it is generally much more
natural, and much easier, to work with formal delta functions, and
with generating functions in general, rather than with individual
components of vertex operators and individual relations among them.
Generating functions of otherwise very complicated objects, such as
nonassociative product operations, or operators on a space, or
identities among such operators, pervade vertex operator algebra
theory and allow one to work efficiently.

There are many generalizations and analogues of (\ref{Jacobi}),
including twisted Jacobi identities, as in \cite{FLM5}; generalized
Jacobi identities for abelian intertwining algebras, etc., as in
\cite{DL1}, \cite{DL2}, \cite{FFR} and \cite{Mos}; Huang's much more
subtle identity for {\it nonabelian} intertwining algebras
\cite{Hua5}; multi-Jacobi identities, as in \cite{Hus1}, \cite{Hus2};
and ``logarithmic analogues'' of the Jacobi identity, both untwisted
\cite{L3} and twisted \cite{DoyLM}, which serve to ``explain'' and
generalize work of S.~Bloch \cite{Blo} on the appearance of certain
values of the Riemann zeta function that arose in certain vertex
operator computations.

We have emphasized that the notion of vertex operator algebra is
actually the ``one-complex-dimensional analogue'' of the notion of Lie
algebra.  But {\it at the same time that it is the
``one-complex-dimensional analogue'' of the notion of Lie algebra, the
notion of vertex operator algebra is also the
``one-complex-dimensional analogue'' of the notion of commutative
associative algebra} (which again is the corresponding
``one-real-dimensional'' notion).  Again, operad language can be used
to make this precise, as we comment below.

{\it The remarkable and paradoxical-sounding fact that the notion of
vertex operator algebra can be, and is, the ``one-complex-dimensional
analogue'' of {\em BOTH} the notion of Lie algebra {\em AND} the
notion of commutative associative algebra lies behind much of the
richness of the whole theory, and of string theory and conformal field
theory.}  When mathematicians realized a long time ago that complex
analysis was qualitatively entirely different from real analysis
(because of the uniqueness of analytic continuation, etc., etc.), a
whole new point of view became possible.  In vertex operator algebra
theory and string theory, there is again a fundamental passage from
``real'' to ``complex,'' this time leading from the concepts of {\it
both} Lie algebra and commutative associative algebra to the concept
of vertex operator algebra and to its theory, and also leading from
point particle theory to string theory.

This analogy with the notion of commutative associative algebra comes
{}from the ``commutativity'' and ``associativity'' properties of the
vertex operators $Y(v,x)$ in a vertex operator algebra, detailed in
\cite{FLM5} and \cite{FHL}, and discussed in many places, including
the book \cite{LL}.  These properties are rooted in
conformal-field-theoretic properties of vertex operators, as in
\cite{BPZ}; see \cite{Go}.  In fact, the Jacobi identity
(\ref{Jacobi}) follows from the commutativity property, in the
presence of certain ``minor'' axioms; see \cite{FLM5}, \cite{FHL},
\cite{Go}.  The term ``commutativity'' actually refers to a certain
``commutativity of left-multiplication operations,'' which is why it
can (and in fact does) imply associativity and the Jacobi identity.

It is natural to ask: {\it Can the Jacobi identity axiom
(\ref{Jacobi}) in the definition of vertex (operator) algebra be
simplified?}  As we have discussed, the Jacobi identity is actually
the generating function of an infinite list of generally
highly-nontrivial identities, {\it and one needs many of these
individual component identities in working with the theory}.  But is
there some ``simpler'' condition that in fact implies the Jacobi
identity (in the presence of the ``minor'' axioms in the definition of
vertex operator algebra)?

In fact there is, and this simpler condition, which is related to the
``commutativity'' property, does indeed look much simpler than the
Jacobi identity axiom, but it turns out that the apparent simplicity
is deceptive.

This simple-looking replacement axiom is:

For all $u,v \in V$, where $V$ is a structure satisfying all the
conditions in Definition \ref{vertexoperatoralgebra} except the Jacobi
identity (\ref{Jacobi}), there exists a nonnegative integer $k$ such
that
\begin{equation}\label{weakcomm}
(x_1 - x_2)^k [Y(u,x_1),Y(v,x_2)]=0.
\end{equation}
This {\it ``weak commutivity''} condition, and also, more
significantly, the theorem that it implies the Jacobi identity (in the
presence of ``minor'' axioms), were discovered in \cite{DL2}, where
(\ref{weakcomm}) was actually a special case of a much more general
condition, namely, the analogous assertion for generalized vertex
algebras and abelian intertwining algebras, mentioned above.  In that
greater generality, formal algebraic functions schematically called
$A$ and $B$ in (\ref{generalizedcommrels}) replace the expression
$(x_1 - x_2)^k$ in (\ref{weakcomm}).  All this is treated in
\cite{DL2}, in the full generality.  The special case (\ref{weakcomm})
is discussed in the Introduction of \cite{DL2}, formula (1.4).  The
proof \cite{DL2} that (\ref{weakcomm}) implies the Jacobi identity,
and that the generalizations of (\ref{weakcomm}) imply the
corresponding generalized Jacobi identities, for generalized vertex
algebras and abelian intertwining algebras, are continuations of the
idea that ``commutativity'' implies the Jacobi identity.  In
\cite{DL2}, we were working with graded structures; just as a vertex
operator algebra is $\mathbb{Z}$-graded, a generalized vertex algebra or
abelian intertwining algebra is graded, too (actually 
$\mathbb{Q}$-graded), and this grading was useful in the proof \cite{DL2} that
weak commutativity and its generalizations imply the corresponding
(generalized) Jacobi identities.  Soon after \cite{DL2}, H.~Li was
able to remove the grading hypothesis, and in particular, he was able
to prove that for a {\it vertex algebra} (without grading), weak
commutativity (\ref{weakcomm}) implies the Jacobi identity, in the
presence of ``minor'' axioms.  This and a number of related results
are covered in \cite{LL}.

Notice that the condition (\ref{weakcomm}) is reminiscent of the
generalized commutation relations discussed above, such as
(\ref{generalizedcommrels}) and in particular, (\ref{crossbracket});
in fact, (\ref{weakcomm}) is of course an example of a generalized
commutation relation.

The fact that such a simple-looking generalized commutation relation
as (\ref{weakcomm}) can serve as an axiom replacing the Jacobi
identity in the definition of vertex operator algebra is not as useful
as it might seem.  In fact, starting in \cite{DL2} itself, we chose
not to take (\ref{weakcomm}) as an ``official'' replacement axiom, in
spite of the fact that we proved, and stated, there that it could be
taken as a replacement axiom.  There are essentially three reasons why
we have chosen not to take (\ref{weakcomm}) as an ``official'' axiom
replacing the Jacobi identity: First, one needs ``all'' the
information in the Jacobi identity (and in the relevant generalized
Jacobi identities, in the context of generalizations of the notion of
vertex operator algebra).  Second, if one wants to prove that a
certain structure $V$ is indeed a vertex (operator) algebra, it is
essentially just as hard to prove the condition (\ref{weakcomm}) as it
is to prove the Jacobi identity, for all $u,v \in V$.  (In other
words, the proof that (\ref{weakcomm}) implies the Jacobi identity is
quite short, and in particular is much simpler than the proof that
either (\ref{weakcomm}) or the Jacobi identity holds for all $u,v \in
V$, for interesting examples of vertex (operator) algebras.)  And
third, (\ref{weakcomm}) fails as a replacement axiom for the notion of
{\it module} for a vertex (operator) algebra.  In general, one can
think of a module for an algebraic structure as a space on which the
algebra acts linearly, such that all the axioms (in the definition of
algebra) that make sense hold; this principle is compatible with the
standard definition of ``module'' for a Lie algebra and the standard
definition of ``module'' for an associative algebra, for instance.
This principle indeed motivates the standard definition of ``module''
for a vertex (operator) algebra, and using (\ref{weakcomm}) in place
of the Jacobi identity would {\it not} lead to the correct notion of
module.  Analogous comments hold in generalized settings such as
abelian intertwining algebras, and also, twisted modules.  For
instance, in \cite{FLM5}, certain twisted Jacobi identities were
proved, and these identities endowed certain ``twisted'' spaces with
what came to be called twisted module structure (although the term
``twisted module'' was not used on \cite{FLM5}); this twisted module
structure was necessary for the construction of $V^{\natural}$.  (I
mentioned above that the vertex operator $Y(x)$ entering into Theorem
\ref{theoremsl2} is an example of a ``twisted vertex operator,'' that
is, a vertex operator appearing in a twisted module.)

As we have been suggesting, in vertex operator algebra theory it is
notoriously difficult to construct nontrival {\it examples} of vertex
operator algebras, even examples that are much simpler than
$V^{\natural}$, and of course one cannot do the theory without
examples; in fact, the theory is so rich because the examples are so
rich.

{\it How can one efficiently construct families of examples of vertex
operator algebras and their modules?}

In classical algebraic subjects like group theory, Lie algebra theory,
etc., one of course has many interesting examples available from the
beginning, guiding the development of the general theory.  In vertex
operator algebra theory, {\it there are no essentially nontrivial
examples that are easy to construct and prove the axioms for.}  (The
only ``easy'' examples of vertex algebras are commutative associative
algebras equipped with derivations; see \cite{B1}.)  This is yet
another reason why vertex operator algebra theory is inherently
``nonclassical.''  In ``classical'' mathematics, there simply {\it
were no nontrivial examples of vertex operator algebras} ``lying
around waiting to be axiomatized,'' in contrast with, say, vector
spaces, groups, Lie algebras, etc., etc.

A conceptually elegant, extremely general and extremely
convenient-to-use solution of the problem of constructing families of
examples of vertex operator algebras, and also modules for them, was
developed by H.~Li (\cite{Li1}, \cite{Li2}), generalizing earlier
constructions of examples of vertex operator algebras, including,
among others, constructions of B.~Feigin-E.~Frenkel \cite{FF} and
I.~Frenkel-Y.~Zhu \cite{FZ} (and constructions in \cite{FLM5}).
Briefly, Li formulated a subtle notion of {\it representation of,} as
opposed to {\it module for,} a vertex operator algebra, by developing
the theory of a ``vertex-algebraic analogue'' of the notion of the
usual endomorphism algebra ${\rm End}\;W$ of a vector space $W$ and by
defining a {\it representation of a vertex algebra $V$} on $W$ to be a
(suitable kind of) homomorphism from $V$ to this endomorphism-algebra
structure on $W$.  This endomorphism-algebra structure has roots in
conformal field theory.  It has also been exploited in work of B.~Lian
and G.~Zuckerman \cite{LZ1}, \cite{LZ2}.  Li's analysis of this point
of view culminated in relatively-easy-to-implement sufficient
conditions for constructing families of vertex operator algebras {\it
and} their modules.  An enhanced treatment of this work of Li's is the
main goal of the book \cite{LL}, which also highlights general
theorems of E.~Frenkel-V.~Kac-A.~Radul-W.~Wang \cite{FKRW},
Meurman-Primc \cite{MP2} and X.~Xu \cite{Xu} useful for constructing
families of examples of vertex (operator) algebras, including those
based on the Virasoro algebra, those based on Heisenberg Lie algebras,
those based on affine Lie algebras, and those based on lattices.  (It
happens, though, that such theorems do not serve to simplify the
construction of the vertex operator algebra $V^{\natural}$, even
though all the kinds of vertex operator algebras just mentioned do
enter into the construction of $V^{\natural}$.)

The analogy between the notion of vertex operator algebra and the
notion of commutative associative algebra is in fact directly related
to the conformal-field-theory viewpoint.  This analogy was pointed out
by I.~Frenkel \cite{F}, in the initiation of a program to construct
(geometric) conformal field theory using vertex operator algebras.  In
\cite{Hua1}, \cite{Hua4}, Huang introduced a precise, and deep,
analytic-geometric notion of ``geometric vertex operator algebra,''
and established that it is equivalent to the (algebraic) notion of
vertex operator algebra.  Formulating and proving the geometric
aspects of the action of the Virasoro algebra was the hard part,
involving differential geometry and analysis on infinite-dimensional
moduli spaces; the infinite-dimensionality comes from the
arbitrariness of analytic local coordinates at punctures on a Riemann
sphere.  The sewing of multipunctured Riemann spheres, with analytic
local coordinates vanishing at the punctures (both the ``incoming''
and ``outgoing'' punctures), is reflected in the structure of
vertex-operator-algebraic operations.  In particular, the {\it formal}
variables are systematically specialized to {\it complex} variables,
and one does extensive analysis in addition to extensive algebra.

An interpretation---actually, restatement---of these constructions and
theorems of Huang, including some discussion of the principle
mentioned above that the notion of vertex operator algebra is a
natural ``complexification'' of the notion of both Lie algebra and of
commutative associative algebra, appears in \cite{HL1}, \cite{HL2}, in
the language of operads.  What is being ``complexified'' is the
(one-real-dimensional) operads underlying both the notion of Lie
algebra and the notion of commutative associative algebra; the new
analytic partial operad underlying the notion of vertex operator
algebra \cite{Hua4} (cf. \cite{HL1}, \cite{HL2}) is
one-complex-dimensional.  In fact, this one-complex-dimensional
partial operad ``dictates'' the algebra of vertex operator algebra
theory.  More precisely, when one takes into account the arbitrary
local coordinates vanishing at the punctures, one really has an
infinite-dimensional structure.  The algebraic operations are
``mediated'' by the infinite-dimensional analytic geometry of this
partial operad, in a precise way; a vertex operator algebra becomes a
{\it representation} of this analytic partial operad.  Again, we are
seeing a mathematical reflection of a passage from point particle
theory to string theory.

In an extensive series of works, K.~Barron (\cite{Ba1}--\cite{Ba6})
has carried out a sophisticated super-geometric
(super-conformal-field-theoretic) analogue of this work of Huang's,
using super-Riemann spheres with general superconformal local
coordinates vanishing at the ``incoming'' and ``outgoing'' punctures,
and using vertex operator superalgebras endowed with the appropriate
super-geometric structure.

As we have mentioned a number of times, in vertex operator algebra
theory, it is crucial to make precise distinctions between formal
variables and complex variables.  This distinction is particularly
dramatic in this work of Huang and Barron.  First they had to carry
out elaborate formal algebra, and then they had to systematically
specialize the (infinitely many) formal variables to complex variables
to obtain the desired results.  In classical mathematics, one is used
to being (appropriately!) ``careless'' about the distinction between
formal and complex variables; for instance, one routinely writes
formal expressions such as $\sum_{n \ge 0}x^n$ without saying much
about whether the variable $x$ is supposed to be formal or complex.
In this simple example, it is so clearly understood that this
geometric series converges for certain complex numbers $x$ and not
others that the notation $\sum_{n \ge 0}x^n$ can easily be used, in
the same discussion, for either the formal sum or the convergent
series, depending on what is being said.  However, in vertex operator
algebra theory, where one also needs both formal and complex
variables, the formal algebra (and ``formal calculus,'' mentioned
above in connection with (\ref{Jacobi}), in which the three formal
variables {\it cannot} all be specialized to complex variables)
becomes so subtle and elaborate that it is necessary to be very
explicit about the distinction between the two kinds of variables.
Correspondingly, in Huang's and Barron's work, there are many formal
theorems, which cannot be initially and directly formulated in terms
of complex variables, and these theorems are then systematically
applied to give analytic consequences.  One particular ``formal''
theorem that was proved by Huang and extended to the superalgebraic
setting by Barron has been considerably generalized in \cite{BHL}, to
a factorization theorem for formal exponentials, in a setting
involving arbitrary infinitie-dimensional Lie algebras.  This theorem
reflects and generalizes the formal algebra required for the sewing of
Riemann spheres or super-Riemann spheres with general coordinates
vanishing at punctures that are sewn together.

Huang's work on geometric vertex operator algebras is a major step of
many in a program to use vertex operator algebra theory and its
representation theory to construct conformal field theories in the
precise sense of G.~Segal's and M.~Kontsevich's definition of the
(mathematical) notion of conformal field theory (see
\cite{S2}--\cite{S4}); this definition was clarified by P.~Hu-I.~Kriz
(\cite{HuK1}, \cite{HuK2}).  Geometric vertex operator algebras amount
to the holomorphic, genus-zero part of the construction of conformal
field theories.  The term ``genus-zero'' now refers to the Riemann
spheres mentioned above; it does not refer to the ``genus-zero
property'' of the discrete subgroups of $SL(2,\mathbb{R})$ discussed
earlier in connection with moonshine.  This use of ``genus zero''
corresponds to conformal field theory at ``tree level,'' in physics
terminology.

A major achievement of this program so far is Huang's proof
(\cite{Hua7}--\cite{Hua11}), in a general setting, of the Verlinde
conjecture and his solution of the problem of constructing modular
tensor categories from the representation theory of vertex operator
algebras.  E.~Verlinde conjectured \cite{Ve} that certain matrices
formed by numbers called the ``fusion rules'' in a ``rational''
conformal field theory are diagonalized by the matrix given by a
certain natural action of the fundamental modular transformation $\tau
\mapsto -1/\tau$ on ``characters.''  A great deal of progress has been
achieved in interpreting and proving Verlinde's (physical) conjecture
and the related ``Verlinde formula'' in mathematical settings, in the
case of the Wess-Zumino-Novikov-Witten models in conformal field
theory, which are based on affine Lie algebras.  On the other hand,
G.~Moore and N.~Seiberg \cite{MSe1}, \cite{MSe2} showed, on a physical
level of rigor, that the {\it general} form of the Verlinde conjecture
is a consequence of the axioms for rational conformal field theories,
thereby providing a conceptual understanding of the conjecture.  In
the process, they formulated a conformal-field-theoretic analogue,
later termed ``modular tensor category'' by I.~Frenkel, of the
classical notion of tensor category for representations of (i.e.,
modules for) a group or a Lie algebra.  A modular tensor category is
in particular a braided tensor category that is also rigid and
``nondegenerate.''

Now, there is a general tensor product theory for modules for a
suitably general vertex operator algebra, a theory based on
intertwining operators \cite{FHL} among modules (the dimensions of
spaces of intertwining operators are the ``fusion rules'' mentioned
above): The tensor product functors and appropriate structure were
constructed in \cite{HL4}--\cite{HL6}, and in \cite{Hua2} Huang proved
a general operator-product-expansion theorem for intertwining
operators, enabling him to construct the natural associativity
isomorphisms between suitable tensor products of triples of modules.
(The existence of such an operator product expansion was a key {\it
assumption}---not {\it theorem}---in \cite{MSe2}.)  The resulting
braided tensor category structure was enhanced in \cite{HL3} to what
we called ``vertex tensor category'' structure.  This structure is
much richer than braided tensor category structure.  Vertex tensor
category structure is ``mediated'' by the analytic partial operad
\cite{Hua4}, \cite{HL1}, \cite{HL2}, based on multipunctured Riemann
spheres with arbitrary analytic local coordinates vanishing at the
punctures, by analogy with how the structure of ordinary, classical,
braided tensor categories is ``mediated'' by operadic structure in one
real dimension.  That is, not only is the concept of vertex operator
algebra {\it itself} ``based'' on the one-complex-dimensional operadic
structure discussed above, but so is a {\it tensor category theory}
\cite{HL3} of modules for a (suitable) general vertex operator
algebra.  In particular, in place of a single tensor product functor,
there is a natural {\it family} of tensor product functors, indexed by
a power of the determinant line bundle over the moduli space of
three-punctured Riemann spheres with analytic local coordinates
vanishing at the punctures, and the natural associativity isomorphisms
among triple tensor products, and the coherence, are controlled by
this geometric structure; this is explained in \cite{HL3}.  In this
tensor product theory, the underlying vector space of the tensor
product of suitable modules for a vertex operator algebra is {\it not}
the tensor product vector space of the modules.  Instead, intertwining
operators among triples of modules form the starting point for a
family of analytically-defined tensor product functors, and it is a
subtle matter to construct the tensor product spaces.

The work \cite{HL3}--\cite{HL6} and \cite{Hua2} was originally
inspired by the work of D.~Kazhdan and G.~Lusztig, starting in
\cite{KLu1}--\cite{KLu3}, constructing a tensor product theory for
certain categories of modules of a fixed non-positive-integral level
for an affine Lie algebra.  However, while the theory of
\cite{HL3}--\cite{HL6} and \cite{Hua2} applies to the module
categories of many families of vertex operator algebras, this theory
does not include \cite{KLu1}--\cite{KLu3} as a special case, because
the modules considered by Kazhdan-Lusztig are not semisimple.  But
recently, in joint work with L.~Zhang \cite{HLZ}, we have generalized
the tensor product theory \cite{HL4}--\cite{HL6}, \cite{Hua2} to
``logarithmic'' tensor product theory, which indeed accommodates
suitable non-semisimple module categories.  It turned out the the work
necessary for this generalization was considerable.  In particular,
the already-intricate formal calculus necessary for
\cite{HL3}--\cite{HL6} and \cite{Hua2} had to be extended to
``logarithmic formal calculus,'' and many of the arguments in
\cite{HL4}--\cite{HL6}, \cite{Hua2} had to be replaced by new ones,
for the generalization.  Using the work \cite{HLZ}, Zhang has
succeeded \cite{Zha} in recovering the braided tensor category of
Kazhdan-Lusztig as a special case of the theory, and of endowing it
with vertex tensor category structure.

To return to Huang's work announced in \cite{Hua7}--- Under very
general and natural conditions on a simple vertex operator algebra
$V$, Huang proved the Verlinde conjecture, and using this result, he
proved the rigidity and in fact modularity of the braided tensor
category constructed in \cite{HL4}--\cite{HL6}, \cite{Hua2}.  Zhu's
theorem \cite{Z}, which I mentioned above, on modular transformation
properties of the graded dimensions of modules, is necessary in the
formulation of the Verlinde conjecture in this general setting, and
Zhu's method of proof plays an important role in the development of
the genus-one theory.  In Huang's proof of the Verlinde conjecture, a
crucial step was to prove the modular invariance of spaces of
multipoint correlation functions involving compositions of
(multivalued) intertwining operators (as opposed to single-valued
vertex operators).  For this, Huang had to develop a new,
analytically-based method \cite{Hua9}, since Zhu's method (mainly, his
use of the commutativity property of vertex operators) cannot be
generalized to this multivalued setting.  (This situation is actually
somewhat analogous to the situation discussed above in connection with
generalized commutation relations (\ref{generalizedcommrels}), where
it was impossible to construct {\it ordinary} commutation relations,
which would have been easier to use, if they had existed.  In Huang's
proof of the Verlinde conjecture, however, the situation is even more
subtle.)  Huang established natural duality and modular invariance
properties for genus-zero and genus-one multipoint correlation
functions constructed from intertwining operators for a vertex
operator algebra satisfying general hypotheses, and as I have
mentioned, the multiple-valuedness of the multipoint correlation
functions led to considerable subtleties that had to be handled
analytically and geometrically, rather than just algebraically.  Then,
with these results having been established, the strategy of Huang's
proof of the Verlinde conjecture reflected the pattern of \cite{MSe1},
\cite{MSe2}; he used these results to establish two formulas that
Moore and Seiberg had derived from strong assumptions, namely, the
axioms for rational conformal field theory, which of course cannot be
assumed here.  It is much harder to (mathematically) prove the axioms
in \cite{S2}--\cite{S4} or the axioms in \cite{MSe2} for conformal
field theory than it is to prove the Verlinde conjecture; in Huang's
work, the truth of Verlinde conjecture was needed to {\it prove} the
desired properties of the tensor category.  In turn, his proof of the
general Verlinde conjecture required a large amount of (existing and
new) representation theory of vertex operator algebras.  The
hypotheses of the theorems entering into Huang's work are very
general, natural and purely algebraic, and have been verified in a
wide range of important examples, while the theory itself is heavily
analytic and geometric as well as algebraic.

In another important direction, a great deal of extensive and deep
work has been done in conformal field theory using algebraic geometry.
I will only mention the books of A.~Beilinson and V.~Drinfeld
\cite{BD}, and of E.~Frenkel and D.~Ben-Zvi \cite{FB}.  The review
\cite{Hua6} of \cite{FB} includes an interesting discussion of the
relation between the algebro-geometric viewpoint and the viewpoint of
the representation theory of vertex operator algebras.

I would like to mention a certain inspiring article by M.~Atiyah
\cite{A}, in which among many other things he compares geometry and
algebra; discusses the interaction between mathematics and physics;
and comments on string theory and also on finite simple groups, and in
particular the Monster and its connections with elliptic modular
functions, theoretical physics and quantum field theory---some of
which we've actually been discussing.  Among his many stimulating
comments are these, about ``the dichotomy between geometry and
algebra'': ``Geometry is, of course, about space... Algebra, on the
other hand...is concerned essentially with time.  Whatever kind of
algebra you are doing, a sequence of operations is performed one after
the other, and `one after the other' means you have got to have
time... Algebra is concerned with manipulation in {\it time}, and
geometry is concerned with {\it space}.  These are two orthogonal
aspects of the world, and they represent two different points of view
in mathematics.  Thus the argument or dialogue between mathematicians
in the past about the relative importance of geometry and algebra
represents something very fundamental.''

In the spirit of what we've been discussing: While a string sweeps out
a two-dimensional (or, as we've been mentioning,
one-complex-dimensional) ``world-sheet'' in space-time, a point
particle of course sweeps out a one-real-dimensional ``world-line'' in
space-time, with time playing the role of the ``one real dimension,''
and this ``one real dimension'' is related in spirit to the ``one real
dimension'' of the classical operads that I've briefly referred to---
the classical operads ``mediating'' the notion of associative algebra
and also the notion of Lie algebra (and indeed, any ``classical''
algebraic notion), and in addition ``mediating'' the classical notion
of braided tensor category.  The ``sequence of operations performed
one after the other'' is related (not perfectly, but at least in
spirit) to the ordering (``time-ordering'') of the real line.  But as
we have emphasized, the ``algebra'' of vertex operator algebra theory
and also of its representation theory (vertex tensor categories, etc.)
is ``mediated'' by an (essentially) one-{\it complex}-dimensional
(analytic partial) operad (or more precisely, as we have mentioned,
the infinite-dimensional analytic structure built on this).  When one
needs to compose vertex operators, or more generally, intertwining
operators, after the formal variables are specialized to complex
variables, one must choose not merely a (time-)ordered sequencing of
them, but instead, a suitable complex number, or more generally, an
analytic local coordinate as well, for each of the vertex operators.
This process, very familiar in string theory and conformal field
theory, is a reflection of how the one-complex-dimensional operadic
structure ``mediates'' the algebraic operations in vertex operator
algebra theory.  Correspondingly, ``algebraic'' operations in this
theory are not instrinsically ``time-ordered''; they are instead
controlled intrinsically by the one-complex-dimensional operadic
structure.  The ``algebra'' becomes intrinsically geometric.
``Time,'' or more precisely, as we discussed above, the
one-real-dimensional world-line, is being replaced by a
one-complex-dimensional world-sheet.  This is the case, too, for the
vertex tensor category structure on suitable module categories.  In
vertex operator algebra theory, ``algebra'' is more concerned with
one-complex-dimensional geometry than with one-real-dimensional time.

As we have discussed, the Monster is indeed (very deeply) related to
string theory.  I mentioned above that one of the distinguishing
features of the moonshine module vertex operator algebra
$V^{\natural}$ is that its representation theory is trivial.  This
means in particular that the braided tensor category attached to it is
also trivial.  But its {\it vertex} tensor category, and
correspondingly, its ``one-complex-dimensional world-sheet'' algebra,
is {\it not} trivial.  With Huang, we have been thinking about whether
it will be possible to realize the Monster geometrically, in terms of
the complex operad.

\vspace{.5in}

\noindent {\small \sc Department of Mathematics, Rutgers University,
Piscataway, NJ 08854} \\
{\em E--mail address}: lepowsky@math.rutgers.edu  
 
\end{document}